\documentclass{article}[12pt]
\usepackage{color,times,amsmath,amsfonts,latexsym,epsfig,epsf,colordvi}
\usepackage{url,hyperref}
\usepackage[english]{babel}

\newcommand{\CC}{\mathbb {C}}
\newcommand{\RR}{\mathbb {R}}

\newcommand{\al}{\alpha }
\newcommand{\bt}{\beta }

\newcommand{\et}{\eta }
\newcommand{\la}{\lambda }

\newcommand{\bp}{\begin{pmat}}
\newcommand{\ep}{\end{pmat}}

\def\tcb{\textcolor[rgb]{0,0,1}}

\author{  Frank Uhlig \thanks{Department of Mathematics and Statistics, Auburn 
University, Auburn, AL 36849-5310 \ (uhligfd@auburn.edu)}}

 \title{\vspace*{-6mm} Zhang Neural Networks for Fast and Accurate  Computations of the Field of Values\\[2mm]}

\begin{document}
\date{~}
\thispagestyle{empty}
\maketitle

\thispagestyle{empty}

\vspace*{-8mm}
\begin{center} { \bf Abstract  } \\[2mm]
\begin{minipage}{150mm}
In this paper a new and different neural network, called  Zhang Neural Network (ZNN) is appropriated from  discrete time-varying matrix problems and  applied  to the angle parameter-varying matrix field of values (FoV) problem. This problem acts as a test bed for newly discovered  convergent 1-step ahead finite difference formulas of high truncation orders. The ZNN method that uses a look-ahead finite difference scheme of  error order 6  gives us 15+ accurate digits of the FoV boundary in record time when applied to hermitean  matrix flows $A(t)$.
\end{minipage}\\[-1mm]
\end{center}  
\thispagestyle{empty}


\noindent{\bf AMS :} 15A60, 65F15, 65F30, 15A18\\[-4mm]

\pagestyle{myheadings}
\thispagestyle{plain}
\markboth{Frank Uhlig}{ZNN for Fast and Accurate Computations of the FoV }

\section*{Introduction }\vspace*{2mm} 

The field of values (FoV)
$$ F(A)\ = \ \{ x^*Ax \mid x \in \CC^n, \|x\|_2 = 1\}  
$$
of real or complex square matrices $A_{n,n}$ has been studied for well over 100 years. In 1902 Ivar Bendixon \cite{Be} gave us the so called Bendixon rectangle that contains the field of values. 
Felix Haussdorf \cite{FH19} and Otto Toeplitz \cite{OT18} independently  proved its convexity exactly 100 years ago. And 40 years ago, Charlie Johnson \cite{J} added rotations to the Bendixon rectangle idea to compute  boundary points of $F(A)$ in the complex plane via repeated eigenvalue computations as accurately as Francis' QR algorithm allows for hermitean matrices, i. e., with 14 to 15  accurate leading digits.
Johnson's method computes the eigenvectors of the extreme eigenvalues of matrix rotations  $A(t)  = e^{i t} \cdot A$ of $A$ via the  hermitean and skew-hermitean parts $H$ and $K$ of $A$, defined as
$$
H \ = \ (A+A^*)/2 \quad \text{ and } \quad  K \ = \ (A - A^*)/(2i) \quad \text{ with } H^* = H \ \text{ and } \ K^* = K \ .
$$
With these notations  any given matrix $A = H + i K \in \CC^{n,n} $ generates the hermitean matrix flow $A(t) = \cos(t)*H + \sin(t)*K $ with $A(t) = A(t)^*$ for all angles $0 \leq t \leq 2\pi$ in its associated FoV problem. The normalized  eigenvectors $x$ for the largest and smallest eigenvalues of $A(t)$ then determine two $\partial F(A)$ points via the quadratic form  $x^*Ax \in \CC$  and two $\partial F(A)$ tangents.  37 years ago Miroslav Fiedler \cite{F} established that the field of values boundary curve $\partial F(A)$ is an algebraic curve except for possible  corner point. Authors such as Rudolf Kippenhahn \cite{RK}, and Marvin Marcus and Claire Pesce \cite{MP} and many others have extended our understanding of the field of values. Workshops on field of values research have been held under the name of 'Numerical Ranges and Numerical Radii' (WONRA)  workshops biannually worldwide since  1992. The most recent WONRA gatherings have attracted researchers in quantum computing theory.\\[-2mm]

This paper deals with recent advances in how to compute the field of values boundary curve $\partial F(A)$ faster and more accurately than Johnson's fundamental eigenvalue method \cite{J} or Loisel and Maxwell \cite{LM} most recent attempt to speed up FoV computations via ODE path following methods.. The author \cite{FU14} previously improved on the accuracy and speed of FoV computations by a geometric optimization process that shrinks the area inside the tangents to the field of values polygon  as determined by the computed Bendixon rectangle  sides and envelops the FoV and the area of the  inner polygon of the computed field of values boundary points that lies inside the FoV due to its convexity. Marcus' and Pesce's idea \cite{MP} of using 2 by 2 compressions of the FoV and their field of values ellipses further 'rounded out' the inner FoV polygon  in \cite{FU14} and led to inner and outer FoV area measurements that converge to the actual field of values  area from above and below upon angle refinements. \\[-2mm]

Matrix eigenvalues and Matrix Theory itself came about from studies of differential operators and self replicating function, called eigenfunctions $f \neq 0$, when studied in the early 1800s. What are the {functions $f$} (if any) for which {${\cal A} f = \alpha f$} for some scalar $\alpha$ and a given differential operator ${\cal A}$? In 1829, Augustin--Louis Cauchy \cite[p. 175]{C1829} applied the eigen concept from differential equation models to matrices and began to view the erstwhile  {`eigenvalue equation'} ${\cal A} f = \alpha f$ as a {`null space equation'}, namely as {${\cal A} f - \alpha f = 0$} or $({\cal A} - \alpha \ id)  f= 0$ for the identity operator with $id\ f = f$ for all $f$. This then led to the use of determinants and trying to find the roots of the characteristic polynomial of $A$. Some of today's newest and fastest FoV computational methods reverse this earlier step of 200 years ago from using differential equations ideas to related concepts for matrices  by again employing differential operators and DE solvers in time-varying matrix problems. Interestingly, they can do so more speedily and more accurately here than standard matrix decomposition  procedures would.\\[1mm]
For example, S\'ebastien Loisel and Peter Maxwell \cite{LM} start with the hermitean part function  $H(e^{it}A) =1/2
 (e^{it}A + e^{-it}A^*)$ for a given $A\in C^{n,n}$.  Clearly $H(e^{it}A)$ is hermitean.
If $\la(t)$ is the largest eigenvalue of $H(e^{it}A)$ and $x(t)$ is its unit
 eigenvector, then $x(t)^*Ax(t)$ lies on $\partial(F(A))$ and this is true  for every $ t \in [0,2\pi]$. Their recent paper \cite{LM} studies the parametrized equation 
\begin{equation} \label{LM1}
H(e^{it}A)x(t) = \la(t)x(t) \ .
\end{equation}
With $S(A)= 1/2 (A - A^*)$ we have $A = H(A) + S(A)$. By requiring all eigenvectors $x(t)$ of $A(t)$ to be unit vectors  and differentiating (\ref{LM1}) with respect to $t$, Loisel and Maxwell \cite[section 6]{LM} obtain the differential equation 
\begin{equation} \label{LM2}
\begin{pmat} x'(t)\\ \la'(t) \end{pmat} = 
\begin{pmat}
H(e^{it}A) - \la(t) I_n & -x(t)\\-x(t)^* & 0 \end{pmat}^{-1} 
\begin{pmat}
-iS(e^{it}A)x(t)\\0 \end{pmat} \ .
\end{equation}
This system of ODEs is treated in \cite[section 6]{LM} as an initial value problem starting from just one explicit set of eigendata computations for $t = 0$ and $H(A)$ itself. The parameter-varying  eigenvalue problem is then solved by following the maximal eigenvalue's path of $H(e^{it}A)$ for $t \in (0,2\pi]$. 
To solve (\ref{LM2}) Loisel and Maxwell use  the Dormand-Prince RK5(4)7M method, see \cite[section 6]{LM}, jointly with Hermite interpolation of degree 4. The chosen integrator has  error terms of order $O(h^5)$ and the method delivers  dense $\partial(F(A))$ output points. \cite{LM} includes detailed descriptions of what might happen if the eigenvalues of $H(e^{it}A)$ cross paths, an occurrence that will happen if $A$ has repeated eigenvalues  and results in piecewise line segment on the field of values boundary curve. In \cite{LM} ways are described to detect such events and how to adjust ODE path following eigen methods in case of eigenvalue crossings. For random entry matrices $A_{n,n} \in \CC^{n,n}$, eigenvalue crossings of the associated matrix flow $A(t)$ do not occur and thus they do not affect  field of values boundary curve computations, see e.g. \cite{NW}. However, normal matrices $A$ have polygonal fields of values and their FoV corners are the matrix  eigenvalues. Normal matrices $A$  can be detected readily from their defining equation $A^*A = AA^*$ and their  FoV  easily determined from their eigenvalues.\\[1mm] 
Here we propose a new, an even faster and more accurate method than \cite{LM} that uses convergent look-ahead finite difference formulas and Zhang Neural Networks, rather than  initial value ODE solvers. For simplicity we restrict our algorithm to non-normal matrices $A$ that do not have  repeated eigenvalues and we do not incorporate  any work-arounds  for matrices with repeated eigenvalues. 

\section{The ZNN Method for Finding Field of Values Boundary Points}

Discretized Zhang Neural Networks are a totally new and different process; they share no properties or steps with matrix decomposition methods nor any with  classical path following or homotopy ideas. Their numerical properties cannot be judged in any form or way by comparisons with known methods for time-varying matrix problems, such as found in the works of Cichocki and Unbehauen \cite{CU92}, Mailybaev \cite{M06}, Dieci and Eirola \cite{DE99}, Wang et al \cite{WCW16} or Yi et al \cite{YFT04} for example.

Any parameter- or time-varying problem that is expressed in compatibly dimensioned matrix and vector equations form can be transformed into a Zhang Neural Network. Assume we want to solve an equation
\begin{equation}\label{ZNNstart}
F(t) = f(A(t),B(t),C(t),...,x(t),y(t), ... ) = g(\al(t),\bt(t),... ) =G(t)
\end{equation}
with time- or parameter-varying functions $A(t),B(t), ...,\al(t), \bt(t),...$ for the unknown $x(t), y(t), ...$.
From the given general equation (\ref{ZNNstart}) and following Cauchy's idea from 1829, we form the homogeneous error equation
\begin{equation}\label{ZNNerror}
E(t) = F(t) - G(t) \stackrel{!}{=} 0
\end{equation}
that ideally should be zero or near zero for all $t$ and the yet unknown time-varying solutions $x(t),y(t),...$. The ZNN Ansatz then requires exponential decay for $E(t)$ by stipulating 
\begin{equation}\label{ZNNOde}
\dot E(t) = - \et E(t)
\end{equation}
with a decay constant $\et > 0$. This translates the given task of solving  (\ref{ZNNstart}) into an ordinary differential equation.\\
 Zhang Neural Networks go back almost 20 years to Yunong Zhang and Jun Wang \cite{ZW01}. These new  methods have been used extensively by engineers and the engineering literature now contains well over 400  papers on problem specific ZNN methods. There are several books on ZNN, such as \cite{ZG2015}. The treated problems are generally time-varying, with parameter inputs $A(t), B(t), ...,\al(t), \bt(t), ...$ coming in at regular frequencies, such as 10, 50 or 100 Hz from sensors on robots or autonomous vehicles and so forth. And their specific tasks may involve optimization, linear solves, operator inversions, eigenvalue computations,  and most any other standard computational set-up that we are familiar with in the static numerical problem solving realm such as Sylvester equations, Lyapunov  equations etc.  Note that static solvers are generally of little use or accuracy in time-varying problems. And besides, it appears that the  numerical theory for parameter-varying (continuous or discrete) numerical problems is not well developed at this time. The new area has not been even understood nor researched from the Numerical Analysis angle of view. Many specific ZNN methods  for time-varying or parameter-varying  problems that have been developed and are successfully used  in applications by engineers today, yet they are still unstudied and mostly unknown in numerical analysis circles, in our textbooks and courses.\\[-2mm]
 
 This paper uses newly developed  high truncation error order convergent look-ahead finite difference schemes from \cite{FU18} and Zhang Neural Networks for the parametrized field of values boundary location problem with 
 \begin{equation}\label{Atdef}
 A(t) = \cos(t) \ \dfrac{A+A^*}{2} + \sin(t) \ \dfrac{A-A^*}{2i} = A(t)^*\in \CC^{n,n} \ \ \  ( A(t)= H(e^{it}A) \text{ in Loisel and Maxwell, \cite{LM}} )\ . 
 \end{equation}
\emph{STEP 1, model making :} For the FoV problem ZNN  starts from the initial eigendata for $A(0) = H$ at $t = 0$ and continues in discrete and predictive time or angle  steps until $t = 2\pi$. For all matrix eigenvalue problems,  Zhang Neural Networks would try to solve the parametric eigenvalue equation 
 \begin{equation}\label{EVequat}
A(t) V(t) = V(t)D(t) \ \text{ for all } \ t \ 
\end{equation}
with two linked time-varying unknowns, the nonsingular eigenvector matrices $V(t)$ and  the diagonal matrices $D(t)$ that contain the eigenvalues of $A(t) \in \CC^{n,n}$. \\
\emph{Step 2; form the error function :} The associated ZNN error equation is 
 \begin{equation}\label{errequat}
E(t) = A(t)V(t) - V(t) D(t)  \stackrel{!}{=} O_n \in \CC^{n,n} \ .
\end{equation}
\emph{Step 2a; specify  the model for one eigenvalue $\la_i$ and eigenvector $x_i(t)$ pair of $A(t)$ :}  Now we consider the error equation eigenvector/eigenvalue pair equation
\begin{equation*}
e_i(t) = A(t)x_i(t) - \la_i(t) x_i(t)  \stackrel{!}{=} o_n\in \CC^n
\end{equation*}
for $ i = 1,...,n$.\\
 \emph{Step 3, stipulate exponential decay for $e_i(t)$ :}  Stipulating exponential decay for $e_i$ means that $ \dot e_i(t) = - \eta e_i(t)$ for each $i$ and -- for symplicity's sake -- a global positive decay constant $\eta$. By expansion and rearrangement this leads to the ODE
\begin{equation}\label{errdeexplii}
 ( A(t) - \la_i(t) I_n) \dot x_i(t) - \dot \la_i(t)  x_i(t)  =   ( - \eta (A(t) - \la_i(t) I_n)   - \dot A(t)) x_i(t) 
 \end{equation}
where $I_n$ is the identity matrix of the same size $n$ by $n $ as $A(t)$. \\
\emph{Step 3a, arrange equation (\ref{errdeexplii}) in augmented matrix form :}  Upon further rearrangement this becomes the matrix ODE
\begin{equation}\label{errdeexplim}
\begin{pmat}
 A(t) - \la_i(t) I_n & -x_i(t)\\
 2 x_i^*(t) & 0 
 \end{pmat}
 \begin{pmat} \dot x_i(t)\\ \dot \la_i(t) \end{pmat}
    =   \begin{pmat} ( - \eta (A(t) - \la_i(t) I_n)   - \dot A(t)) x_i(t) \\ -\mu (x_i^*(t)x_i(t) -1) \end{pmat} \ ,
 \end{equation}
where the bottom block-matrix row ensures that the computed eigenvectors become unit vectors with stipulated exponential error decay for a separate decay constant $\mu$. \\  
These are the three ubiquitous paper and pencil start-up steps for all  ZNN based matrixmethods.\\[-3mm]

In our FoV computations, the matrix flow $A(t)$ consists entirely of hermitean matrices and therefore all occurring eigenvalues $\la_i(t)$ will be real and the top left block $A(t) - \la_i(t) I_n$ of the system matrix in (\ref{errdeexplim}) is hermitean for all $t$. By dividing the lower block row in (\ref{errdeexplim}) by --2, we obtain the following  equivalent model, now with a completely hermitean system matrix for speedier computations:
\begin{equation}\label{errdeexplimherm}
\begin{pmat}
 A(t) - \la_i(t) I_n & -x_i(t)\\
 - x_i^*(t) & 0 
 \end{pmat}
 \begin{pmat} \dot x_i(t)\\ \dot \la_i(t) \end{pmat}
    =   \begin{pmat} ( - \eta (A(t) - \la_i(t) I_n)   - \dot A(t)) x_i(t) \\ \mu/2 \cdot (x_i^*(t)x_i(t) -1) \end{pmat} \ .
 \end{equation}
Further details of this process are contained in \cite[section 2]{YZMeZNN}.\\
 The  system matrix $P(t_k)$ above is identical to the one of equation (\ref{LM2}) that was obtained by  differentiating the matrix eigenequation (\ref{LM1}) and was used by Loisel and Maxwell's  \cite{LM} IVP ODE path following method. But the right hand side  of our two coupled matrix ODEs in (\ref{errdeexplimherm}) is  completely different and besides, our  ODE version (\ref{errdeexplimherm}) with a hermitean system matrix contains the decay constants $\et$ and $\mu$ and completely different  terms than (\ref{LM2}).\\
The continuous-time matrix eigendata model (\ref{errdeexplim}) has been proven convergent and robust for noisy inputs in \cite{ZYLUH}  for symmetric matrix flows $A(t) = A(t)^T$. The ideas and proofs of \cite{ZYLUH} easily generalize to hermitean flows $A(t) = A(t)^*$.\\
 Here we deal with the specific discretized version of the eigendata FoV model that we have just introduced. ZNN methods  for more general time-varying matrix problems were studied and tested earlier in \cite{YZMeZNN}.\\[-4mm]
 
  For discretized input data matrices $A(t_k) \in \CC^{n,n}$ and $k = 0,1,2,...$ it is most natural to discretize their differential equations  in sync with the sampling gap $\tau = t_{k+1} - t_k$.  In the following we assume  that $\tau$ is constant  for all $k$.   With
\begin{equation}\label{Pzqdef}
\begin{array}{c}
 P(t_k) = \begin{pmat}
 A(t_k) - \la_i(t_k) I_n & -x_i(t_k)\\
 - x_i^*(t_k) & 0 
 \end{pmat} \in \CC^{n+1,n+1} , \  \ \  z(t_k) = \begin{pmat}  x_i(t_k)\\  \la_i(t_k) \end{pmat} \in \CC^{n+1} \ , \\[3mm]
 \text{ and }  \ q(t_k) = \begin{pmat} ( - \eta (A(t_k) - \la_i(t_k) I_n)   - \dot A(t_k)) x_i(t_k) \\ \mu/2 \cdot (x_i^*(t_k)x_i(t_k) -1) \end{pmat}  \in \CC^{n+1}
 \end{array}
 \end{equation}
 our model (\ref{errdeexplimherm}) with a hermitean system matrix flow $P(t)$ becomes  the matrix ODE 
 \begin{equation}\label{Pzqde}
 P(t_k) \dot z(t_k) = q(t_k) 
 \end{equation}
for $k = 0,1,2,...$, each equidistant discrete time step $0 \leq t_k \leq t_f$ and for each associated eigenpair $x_i(t_k)$ and $\la_i(t_k)$ of $A(t_k)$. Here $t_f$ denotes the final time $t_f$.\\[1mm] 
Following the above start-up steps,  we are now ready to work through the discretized ZNN finite difference portion of the scheme. The ZNN method ultimately relies entirely on look-ahead convergent finite difference schemes. The ensuing high accuracy predictions of future states of a system differ from any other method that we have ever seen or encountered and are completely new. Therefore ZNN methods behave quite differently for time-varying matrix problems. They belong to a special, a new branch of numerical analysis that deserves to be studied without prejudice in its own right.\\[1mm]
To finish solving (\ref{Pzqde}) for  matrix  problems, discrete ZNN methods rely on high error order and convergent  look-ahead finite difference methods. These neural networks generally show great benefits for time-varying problems. Alternately, one might try to use static methods and compute the eigendata of $A(t_k)$ at each time instant $t_k$. Thereby one would know the system properties at time $t_k$ shortly after time $t_k$ has passed due to the inherent computational lag. For modern robot control and autonomous vehicle  and in other applications, however, such knowledge is of little value and might even be dangerously deceptive. When we walk or drive along a path or street it is good to know where we currently are, i. e., our static data or where and  how we have just been. But much more importantly, when moving about we  need to look ahead, even if only  somewhat unconsciously,   and  adjust our direction, our speed and so forth  at time $t_k$ for the future time instance $t_{k+1}$. To do so at $t_k$, we can only rely on the known system data from past or current moments $t_j$ with $j \leq k$. As we progress on our path we must stay on the road, follow its curves, and try to avoid puddles and potholes etc that may lie ahead. Thus we need to make accurate assessments of what comes next upon us in the immediate future using our present and  past data. Clearly meaningful numerical methods for such time-varying problems  must  be able to predict a  system's states in the future at time $t_{k+1}$ and do so reliably from data for times at or before $t_k$. Discrete Zhang Neural Networks rely on convergent look-ahead finite difference schemes that can predict the system's subsequent states reliably and accurately at time $t_{k+1}$ from past data. ZNN methods can do so quickly  even for corrupted data, corrupted by sensor or transmission failures, thanks to their stipulated exponential error decay; for examples see \cite{YZMeZNN}. 
Moreover the convergence rate and accuracy of discretized ZNN methods depend solely on the chosen sampling gap $\tau$, on the chosen decay constants $\et$, $\mu$, ... , and on the forward difference formula's truncation error order.  Backward finite difference formulas may also be employed inside ZNN methods, such as, for example,  to approximate the derivative $\dot A(t_k)$ in the right hand side computations of $q(t_k)$ in formulas  (\ref{Pzqde}) and (\ref{Pzqdef}) above.\\

Look-ahead finite difference formulas are rare in the literature and in our handbook lists such as in \cite[section 14.2]{EU96}. Moreover the few that are listed and involve states at times $t_{k+1}, t_k, t_{k-1}$ through $t_{k-j}$ for some $j \geq 1$ are usually not convergent and thus of no help with discretized ZNN methods that excel at predicting or computing future states fast and accurately. Convergence of multistep formulas is standardly defined via their characteristic polynomials. A characteristic polynomial  is convergent if it is zero stable, if all its characteristic polynomial  roots lie in the closed unit disk of $\CC$, and  if all roots on the unit circle are
 simple, see e.g. \cite[section 17.4]{EU96}. In the past,  only six {\em convergent} look-ahead finite difference schemes have been known, all with relatively low truncation error orders of at most  four, see e.g. \cite{FU18}.
These were either of classical origin, such as Euler's formula or constructed by hand from Taylor expansions for  the iterates $z_{k+1}, z_{k-1}, ..., z_{k-\ell}$ with clever manipulations  that would cancel  second and higher order derivatives in a linear combination of these Taylor expansion equations. For details we refer to \cite{FU18}. Recently this author \cite{FU18} has constructed an algorithm for finding convergent look-ahead finite difference formulas of arbitrary truncation error orders.   New such finite difference schemes have been found up to truncation error order eight. In this paper we will apply and test several convergent look-ahead finite difference formulas to the matrix FoV boundary problem; one known of truncation error order 4, found by hand, pencil and paper in \cite{LMUZ}, another new one of order 6 from \cite{FU18}, as well as four more such formulas of truncation error orders 3, 5, 7 and 7, respectively. Recall that this paper  and the FoV matrix problem is intended as a test bed for comparisons of fast and accurate ZNN based computations. Our ZNN methods achieve spectacular results, see the next section. \\[-6mm]

\section{Field of Values Computations via ZNN and Look-Ahead Finite Difference Formulas}

\emph{ Step 4, clear up all variables appearing  in or after Step 3a :} The iterates of our algorithm are the entries in\\[-2mm] 
$$z_k = z(t_k) = \begin{pmat}  x_i(t_k)\\  \la_i(t_k) \end{pmat} \in \CC^{n+1} $$
from (\ref{Pzqdef}). Here $x_i(t_k)$ denotes a field of values generating extreme eigenvector of $A(t_k) = \cos(t_k)H + \sin(t_k)K$ as defined in (\ref{Atdef}) and   $\la_i(t_k)$ is the corresponding eigenvalue of $A(t_k)$ for the  angle $t_k \in [0,2\pi]$. 
For non sensor based processes such as our formulaic equations based FoV boundary point problem, the derivative $\dot A(t_k)$ that occurs inside the ZNN  based ODE (\ref{Pzqde}) in $q$ and formula (\ref{Pzqdef})  does not need to  be approximated by a backward finite difference formula. 
The derivative of  $A(t_k) = \cos(t_k) H +  \sin(t_k) K$  with respect to $t$ is simply $\dot A(t_k) = \sin(t_k)H - \cos(t_k)K$.\\[1mm]
\emph{Step 5, computer implementation of a difference scheme with fictitious name }\verb^m_s^ :  
Here $m \geq s$ are both positive integers. The new eigendata iterate $z_{k+1} = z(t_k+\tau)$ at time $t_{k+1}$ is  found in ZNN from $\dot A(t_k)$,  from previous eigendata for $t_j$ with $j \leq k$, the respective previous eigenvectors {\tt zs}, and the associated previous eigenvalues {\tt ze}  that are stored in {\tt ZZ} $\in \CC^{n+1,m+s}$, as well as  the $n+1$ by $n+1$ matrix $P$ as detailed in (\ref{Pzqdef}).  This is expressed in the following MATLAB code lines:\\[-5mm]

{\small
 \begin{verbatim}
 % top eigenvalue iteration :
   ze = ZZ(n+1,1);
   zs = ZZ(1:n,1);
   Al = Ath; Al(logicIn) = diag(Al) - ze*ones(n,1); % Al = Ath - ze*In
   P = [Al,-zs;-zs',0];                             %  P is hermitean
 \end{verbatim}}
  
  \vspace*{-5mm}
  \noindent
 The  two code lines below define the right hand side vector $q(t_k)$ in (\ref{Pzqdef}) and find the solution $X$ with $P(t_k) X = q(t_k)$ as indicated in equation (\ref{Pzqde}):

{\small
 \begin{verbatim}
  q = taucoeff*tau*[(eta*Al + Adot)*zs; -1.5*eta*(zs'*zs-1)]; % n+1 vector q in (11)  
  X = linsolve(P,q);  % hermitean matrix P in linsolve without any Lopts.SYM setting
\end{verbatim}}

\vspace*{-1mm}
\noindent
Following some experimentations, we use the increased decay constant $\mu =3 \et$ in the bottom entry of $q$ to achieve  faster  convergence of the eigenvectors to unit length. Finally the eigendata vector $z_{k+1}$ is evaluated as {\tt Znew} in the next code line via the chosen finite difference formula's characteristic polynomial of degree $m+s$, the solution vector {\tt X} of $Px = q$  and the previous eigendata matrix {\tt ZZ} for $z_j$ from $m+j$ indices $j \leq k$, namely \\[1mm]
\verb^  Znew = -(X + ZZ*polyrest);         % from the formula for zdot ^.
\\[1mm]
The only change from one ZNN method to another in this code involves the choice of which normalized finite difference polynomial to use. The respective  characteristic  polynomial is stored -- without its leading coefficient 1 -- in \ \verb^polyrest^ $\in \RR^{m+s}$. \\[1mm]
\emph{Step 6, choosing a convergent 1-step ahead finite discretization formula :} Below we test our ZNN method for two different convergent polynomials and their associated look-ahead finite difference formulas, namely the \verb^2_2b^ one with truncation error order 4 and the \verb^4_5a^ formula with truncation error order 6. The polynomial   \verb^2_2b^  was found by pencil and paper. It has integer coefficients  and first appeared in \cite{LMUZ} while  \verb^4_5a^ was  computed more recently via the MATLAB optimization and search code detailed in \cite{FU18}. Both schemes work well for  FoV boundary computations. The difference scheme \verb^2_2b^ has the non-normalized characteristic polynomial  $p_4(x) = 8x^4+ x_3-6x^2-5+2$ and the discrete  ZNN method \verb^2_2b^ with $p_4$ needs $m+s = 2 + 2 = 4 $ starting values of both the eigenvectors and eigenvalues for $A(t_k) $ at $k = 0, ..., 3$. These we compute at start-up using explicit Francis QR eigenanalyses. The roots of $p_4$ are  $-0.71597 \pm 0.54945i,  1.0000,$ and $ 0.30693$ which lie well within the unit disk in $\CC$, see Figure 2. Our computer generated finite difference scheme \verb^4_5a^ has a non-rational  characteristic polynomial which we approximate in MATLAB double precision as

\vspace*{-5mm}
{\small
\begin{eqnarray*}
p_9(x) &=&-1.632891580619644x^9 -1.084874852377588x^8 + 1.514338299609167x^7 +2.121238162639099x^6 \\
&&-0.3010929138446914x^5-0.9393487657815317x^4 +0.06714730122560907x^3\\ &&+0.3319027505915695x^2 
 -0.04244088319409350x -0.03397751824789656 \ . 
\end{eqnarray*}}

\vspace*{-5mm}
\noindent
Method \verb^4_5a^  uses $p_9$ in its normalized form and it needs $4+5 = 9$ eigendata starting values, which are again found via Francis QR. The roots of $p_9$ are $-0.65551 \pm 0.62479i,
-0.70266 \pm 0.34906i,  1.0000,$ $ 0.43671 \pm 0.32838i,0.47141,$ and $  -0.29289 $. These all lie properly inside the unit disk in Figure 2, insuring convergence when used inside a look-ahead finite difference based ZNN scheme.\\[-2mm]

\noindent
\emph{Step 7, run tests :}
All our test tables for  ZNN field of values computations below  use random entry complex matrices $A_{n,n}$ of varying dimensions. In our first two test tables the matrix dimensions  $n$  increase in size by  factors of 3  and all tables use  sampling gaps  such as $\tau = 0.01 (100 H\!z), 0.005, 0.001, 0.00015$ or $0.0001$. The most accurate  field of values boundary points that anyone can create in MATLAB  -- without use of accuracy enhancements such as inverse iteration or extended digital precision ranges -- are those from  Francis QR algorithm eigenanalyses of  $A(t_k)$. These are accurate in their leading 14 or 15 digits. Loisel and Maxwell \cite{LM} interpolate their ODE computed results further. We instead take smaller and smaller sampling gaps $\tau$ for the ZNN iterations and generate up to 62,833 FoV boundary points via ZNN when $\tau = 0.0001$. This many points seem to be enough for an accurate  polygon approximation of the FoV boundary. The accuracy  is measured in our code by the average  number of   correct  ZNN computed digits of  the boundary point estimates $x^*Ax \in \CC$  for $x = (z_{k})_ {[1:n]}\in \CC^n$,  the eigenvector part of  our iterate $z_{k}$ that was computed predictively at time $t_{k-1}$ via ZNN, when compared  with the Francis QR computed FoV boundary point for $A(t_{k})$.\\
 Note again that this predictive way of computing by design with ZNN is not available in any other time-varying matrix algorithm.\\[-3mm]

\vspace*{2mm}
\enlargethispage{10mm}
{\small 
\hspace*{-6.5mm}
\begin{tabular}{|c||c|c|c|c|c|} \hline &&&&\\[-3mm]
{\bf\verb+ 2_2b + finite}&$n = 3$&$n=9$&$n=27$&$n=81$&$n=243$\\
{\bf diff. formula}&$h \in [0.085,0.093]$&$h \in  [0.085,0.09]$&
$h \in [0.065,0.09]$&$h \in [0.035,0.09]$& $h \in [0.035,0.09]$\\[1mm]
\hline  & & & & &\\[-3mm] \hline & & & & &\\[-3mm]
&$\et = 19$&$  \et = 19$&$ \et = 13$&$\et = 7$&$\et =7$\\
$\tau = 0.005$& 8.5 acc. digits&  7.9 acc. digits& 7.7 acc. digits& 7.2 acc. digits & 7.2 acc. digits\\
1259 FoV points&0.06 sec&0.13 sec&0.15 sec&0.55 sec&3.68 sec\\[1mm] \hline &&&&  \\[-3mm]
&$\et = 80$&$  \et = 90$&$ \et = 84$&$\et = 89$&$\et =90$\\
$\tau = 0.001$& 11.2 acc. digits&  10.7 acc. digits& 10.6 acc. digits& 10.3 acc. digits & 10.4 acc. digits\\
6283 FoV points&0.25 sec&0.44 sec&0.57 sec&2.5 sec&17.3 sec\\[1mm] \hline &&&&  \\[-3mm]
&$\et = 930$&$  \et = 900$&$ \et = 900$&$\et = 900$&$\et =900$\\
$\tau = 0.0001$& 14.9 acc digits& 14.6 acc. digits&  14.3 acc. digits& 14.0 acc. digits& 14.1 acc. digits\\
62833 FoV points&2.1 sec&2.6 sec&5.2 sec&24.3 sec&171 sec\\[1mm] \hline   \\[-3mm]
 & \multicolumn{4}{c|}{ }& $\et = 600$\\
 $\tau = 0.00015$& \multicolumn{4}{c|}{ }& 13.5 acc. digits\\
41889 FoV points& \multicolumn{4}{c|}{ }&  114 sec\\[1mm]
\hline 
\end{tabular}\\[2mm] }
Table 1 : Optimal $\et$ values, accurate digits and run times for the ZNN method   \verb^2_2b^ of truncation error order 4.\\[-2mm]

Next the analogous  data table for ZNN and the field of values boundary problem when using the finite difference formula   \verb^4_5a^   of truncation error order 6.\\[3mm]
\vspace*{2mm}
\hspace*{1mm}
{\footnotesize
\begin{tabular}{|c||c|c|c|c|c|} \hline &&&&\\[-3mm]
{\bf\verb+4_5a + finite}&$n = 3$&$n=9$&$n=27$&$n=81$&$n=243$\\
{\bf diff. formula}&$h \in [0.031,0.041]$&$ h \in [0.037,0.05]$&
$h \in [0.042,0.053]$&$h \in [0.0445,0.051]$& $h \in [0.042,0.057]$\\[1mm]
\hline  & & & & &\\[-3mm] \hline & & & & &\\[-3mm]
&$\et = 6.2$&$  \et = 7.4$&$ \et = 8.3$&$\et = 8.9$&$\et =8.4$\\
$\tau = 0.005$& 10.8 acc. digits&  8.5 acc. digits& 9.5 acc. digits& 9.3 acc. digits & 9.0 acc. digits\\
1259 FoV points&0.05 sec&0.07 sec&0.15 sec&0.58 sec&3.7 sec\\[1mm] \hline &&&&  \\[-3mm]
&$\et = 41$&$  \et = 49$&$ \et = 47$&$\et = 43$&$\et =49$\\
$\tau = 0.001$& 14.7 acc. digits&  13.7 acc. digits& 13.3 acc. digits& 13.0 acc. digits & 13.2 acc. digits\\
6283 FoV points&0.21 sec&0.25 sec&0.59 sec&2.7 sec&16.6 sec\\[1mm] \hline &&&&  \\[-3mm]
&$\et = 390$&$  \et = 500$&$ \et = 530$&$\et = 510$&$\et =520$\\
$\tau = 0.0001$& 15.4 acc digits& 15.3 acc. digits&  15.3 acc. digits& 15.3 acc. digits& 15.2 acc. digits\\
62833 FoV points&1.87 sec&2.3  sec&5.6 sec&26.7 sec&164 sec\\[1mm] \hline   \\[-3mm]
 & $\et = 290$ &$\et = 360$&$\et = 380$&$\et = 370$& $\et = 380$\\
 $\tau = 0.00015$& 15.3 acc. digits&15.3 acc. digits& 15.3 acc. digits&15.1 acc. digits& 15.3 acc. digits\\
41889 FoV points& 1.26 sec &1.55 sec& 3.7 sec&18.0 sec&  111 sec\\[1mm]
\hline 
\end{tabular}\\ }
Table 2 : Optimal $\et$ values, accurate digits and run times for ZNN method   \verb^4_5a^ of truncation error order 6.\\[-3mm]

To round out our comparisons, we next look at four convergent look-ahead finite difference schemes of different error orders and the random entry complex matrix FoV boundary problem.

\vspace*{3mm}
\hspace*{-6mm}
{\small
\begin{tabular}{|c||c|c|c|c|} \hline &&&&\\[-3mm]
{\bf Complex random entry}&\verb+1_2+ f. diff. formula& \verb+3_3+ f. diff. formula&\verb+5_5+ f. diff. formula& \verb+5_6b+ f. diff. formula\\
{\bf matrix $A_{n,n}$ with $n = 243$}&$0.22 \leq h \leq0.26$&$ 0.027 \leq h  \leq0.032$&
$ h = 0.0069$& $h = 0.016$\\[1mm]
\hline  & \\[-3mm] \hline  &  \\[-3mm]
 $\tau = 0.01$ & $\et = 22$ &\multicolumn{3}{c|}{ }\\
(100 Hz)& 5.6 acc. digits &\multicolumn{3}{c|}{ }\\
 631 FoV points & 1.9 sec &\multicolumn{3}{c|}{ }\\[1mm]
\hline \\[-3mm]
 &$\et = 260$&$  \et = 27$&$ \et = 6.9$&$ \et = 16$\\
$\tau = 0.001$& 8.7 acc. digits&  11.6 acc. digits& 11.5 acc. digits& 13.7 acc. digits\\
6283 FoV points&16.6 sec&16.5 sec&17.0 sec& 16.6 sec\\[1mm] \hline &&&  \\[-3mm]
&$\et = 2500$&$  \et = 320$&$ \et = 69$&$ \et = 160$\\
$\tau = 0.0001$& 11.7 acc. digits&  15.0 acc. digits& 14.7 acc. digits& 15.2 acc. digits\\
62833 FoV points&165 sec&166 sec&169 sec& 167 sec\\[1mm] \hline   
\end{tabular}\\[3mm] }
Table 3 : Optimal $\et$ values, accurate digits and run times for the ZNN methods  
 \verb^1_2^ ,  \verb^3_3^ , \verb^5_5^  and  \verb^5_6b^ of \linebreak \hspace*{13.3mm} 
 truncation error orders  3, 5, and 7 twice, respectively.\\[-3mm]
 
 The row data in Table 2 shows that  the ZNN methods  based on the polynomial  \verb^4_5a^   with $\tau = 0.00015$ or $\tau = 0.0001$ give us more than 15 accurate leading digits for 40,000+ or 60,000+ FoV boundary points, respectively, for all chosen matrix sizes $n$ by $n$.  CPU times increase by factors between 1.2 and 7 times for each tripling of  test matrix dimension $n$ from 3 to 243 and for each fixed $\tau$. For any fixed $\tau$ the number of accurately computed digits remains almost the same for all tested dimensions, as does their optimal $\et$ decay constant at around 300 and 500, respectively.  Looking at the column data in Table 2, i. e. for decreasing sampling gaps $\tau$ and  fixed matrix dimension $n$, for  \verb^4_5a^  the auxiliary constant $h = \et \cdot \tau$ is bounded by 0.03 and 0.06 in each column. \\[1mm]
 The data for ZNN with   \verb^2_2b^ in Table 1 shows  similar qualities as that of Table 2, albeit with different optimal $\et$ progressions and generally lesser   numbers of accurately computed digits (by 1 or 2 digits) in their computed FoV boundary points. Table 3 explores the behavior of four adjacent  methods with  truncation error orders $O(\tau^3), O(\tau^5)$ and $O(\tau^7)$ (twice) and for dimension $n = 243$ only. The \verb+3_3+ based ZNN method with truncation error order $O(\tau^5)$ behaves nearly identically to \verb^4_5a^ but still comes short of \verb^4_5a^'s better accuracy performance. Method \verb+3_3+ runs optimally with a much smaller $h = 0.0069$ than that of \verb^4_5a^. And finally, the ZNN methods with the convergent look-ahead finite difference formulas \verb^5_5^ and \verb^5_6b^, both of truncation error orders $O(\tau^7)$  do not  improve significantly on the results of \verb^4_5a^ or \verb^3_3^ despite their respective higher truncation error order. Note however, that \verb^5_6b^ works much better regarding accurate digits than \verb^5_5^ does and compare the different root location geometry  of their associated characteristic polynomials in Figure 2 below.\\[1mm]
Further note that  our lowest truncation order method  \verb^1_2^ is the only one that allows
 usable output with $100H\!z$ discrete data input. The other tested ZNN methods would not converge with $\tau$ much below $0.001$ or $1000H\!z$. Also, the value of $h = \et \cdot \tau$  is well below 0.1 for all methods except for \verb^1_2^ with $h \approx 0.25$ for achieving optimal data accuracy in  FoV  applications. Usually the values of $h$ hover around 1 to 10 in time-varying ZNN methods. The first 'recurrent  neural network' paper with ZNN exponential error decay stipulations by Yunong Zhang and Jun Wang \cite{ZW01} includes the following explanation for the 'small $h$' phenomenon that we have observed in our FoV application. In \cite[overlapping sentence on p. 1163, 1164]{ZW01} there is the following
 remark on  combined networks, such as ours with two independent error equations. ``The
( ... ) process using neural networks is virtually a multiple time-scale 
  system. The plant and controller are in a time scale. $N_z$ and $N_k$ are in another time scale smaller than that of the plant and controller, which requires that the dynamic process for computing KH be sufficiently quick, or the  plant be sufficiently slow time-varying." Thus our aim to achieve high global accuracy of the ZNN computed FoV boundary points by using very small sampling gaps throughout the algorithm counteracts the problem specifics for our chosen knife blade shaped fields of values with long, rather straight sides and very tight turns at the straights' ends for large $n$. The insight that Zhang neural networks are seemingly challenged with our necessarily quick, narrowly spaced angle changes in high FoV curvature areas and slow angle-variances across the straights helps explain our small and very sensitive $h$ anomaly due to our specific desires here. Yet ZNN still comes through very well.\\[-4mm]
  
 Loisel and Maxwell \cite{LM} stop their ODE solver based path following method at 14.5 digit accuracy. This takes  around 405 seconds of CPU time for a 250 by 250 random entry test matrix $A$, see \cite[Fig. 8.1,  +++ curve]{LM}. Our \verb^4_5a^ ZNN method delivers around 15.3 digits accurately in 111 seconds for  243 by 243 random entry matrices when it computes 41,889 FoV boundary points for $\tau = 0.00015$. Here is the MATLAB  loglog line plot for Loisel and Maxwell's path-following ODE method \cite[Fig. 8.1]{LM} and our ZNN based  method \verb^4_5a^. Figure 1 depicts run times versus algorithm accuracy and includes the accuracy/CPU time point \tcb{*} for Johnson's method \cite{J} when used with MATLAB's recently optimized Francis {\tt eig} m file with multishifts.\\[-6mm]
 \begin{center}
\includegraphics[width=98mm]{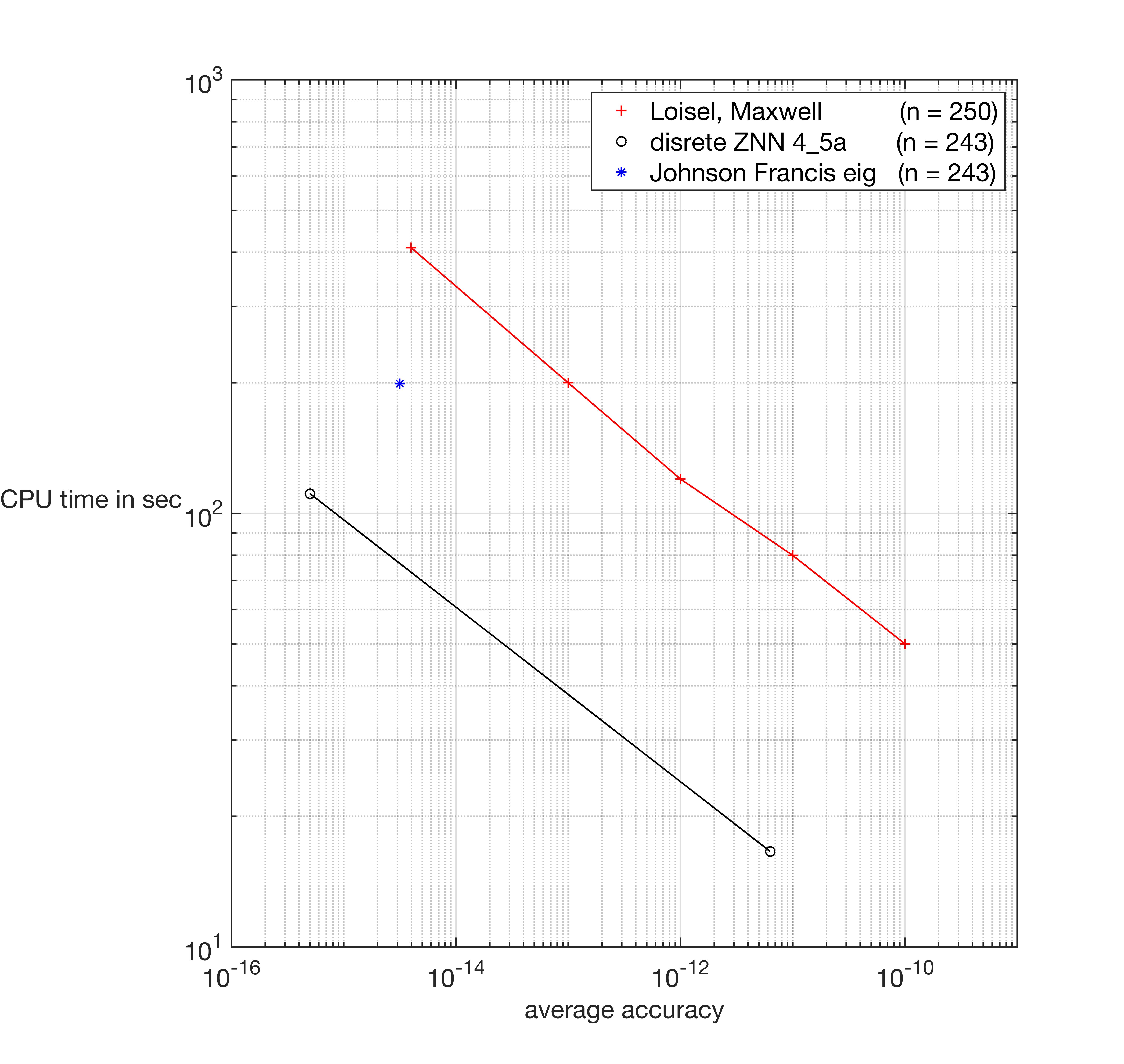} \\[-4mm]
Figure 1
\end{center} 

\vspace*{-2mm}
Figure 1 shows that our  discrete ZNN method with  \verb^4_5a^ is about 4 times faster at each accuracy level than the current best FoV problem solver, i.e., the ODE path following method of \cite{LM}. To draw the  ZNN accuracy data line in Figure 1 we have used the \verb^4_5a^ run time data for $\tau = 0.001$ and $n = 243$ from Table 2 with a 13.2 digit accuracy achieved in 16.6 seconds from 6,283 computed FoV boundary points. Here we have allowed for 2 accurate digits lost  -- from the 13.2 actual accurate digits  on average at 6,000+ FoV boundary points down to just 11.2 to draw the ZNN line in  Figure 1.  This accounts for the unavoidable errors in the gaps between the computed  boundary points due to the natural bulging out of the true convex boundary FoV curve. \\[-4mm]

  It is interesting to note here that specialized  MATLAB programs that rely on  particular eigendata acquisition methods for special matrix problems  seem to have become  almost  irrelevant now due to  recent multishift improvements of {\tt eig} in MATLAB. Using the built-in MATLAB m file {\tt eig} today in Johnson's original FoV algorithm \cite{J} to acquire 42,000 FoV boundary points clearly bests   Loisel and Maxwell \cite{LM} in speed by a factor of two  and achieves the same level of accuracy. In turn today's  {\tt eig} MATLAB function is  bested by discretized ZNN, again  by a  factor of around two according to the CPU data and  Figure 1. This paper was intended to test, examine and show differences in our newly discovered high order discretization versions of ZNN. What further speed gains are possible with local on-chip implementations of discretized ZNN methods now begs an answer.\\[-4mm]
 
 From our set of data for various error order ZNN methods, it is obvious that  ZNN methods of lesser error orders can  give us lower accuracies, but their CPU times are all  similar to those of the higher order methods \verb^4_5a^ or \verb^5_6b^. This is due to the low  $O(n)$ CPU work of ZNN to generate future FoV boundary points via finite difference schemes as  linear  combinations of previous vector data. The length of an $n$ vector does not really matter much for  $O(n)$ computational processes since in our ZNN algorithms, the main effort and CPU time is spent   on solving the linear systems (\ref{Pzqde}) $ P(t_k) \dot z(t_k) = q(t_k)$  in the MATLAB code line \verb^ X = linsolve(P,q)^  that was quoted at the start of this section. This line, when executed 40,000+ times for $n = 243$ with $\tau = 0.00015$, takes up almost half (44 \% according to MATLAB's profile viewer) of the total CPU runtime for \verb^4_5a^. How to solve linear equations more efficiently, possibly via ZNN methods, than using $O(n^3)$ Gauss or generalized Cholesky or $O(n^{2.8...)}$ Strassen inspired methods is an intriguing open problem.\\[-3mm]

The author's MATLAB function {\tt runconv1step1Plot(runs,jend,k,s)} in \cite{FU18} and \cite{FUconv1step1Plot} creates convergent normalized polynomials for possible use in look-ahead finite difference formulas.  It does so from a seed vector $y$ of length $s\geq k$ that the operator chooses and then it creates such polynomials, named \verb^k_s^ here, with truncation error orders $k+2$ for the associated difference scheme. Better high order finite difference methods than ours are likely to exist. How to find or create such is unknown. We feel that the quality of a  polynomial in finite look-ahead difference methods might depend on the lay of the polynomial's roots inside the unit disk.\\[1mm]
To help contemplate these matters, here are the root plots of the six characteristic  polynomials for our the convergent look-ahead discretizations that we have chosen to create Tables 1, 2, and 3.\\[-6mm] 

 \begin{center}
\includegraphics[width=79mm]{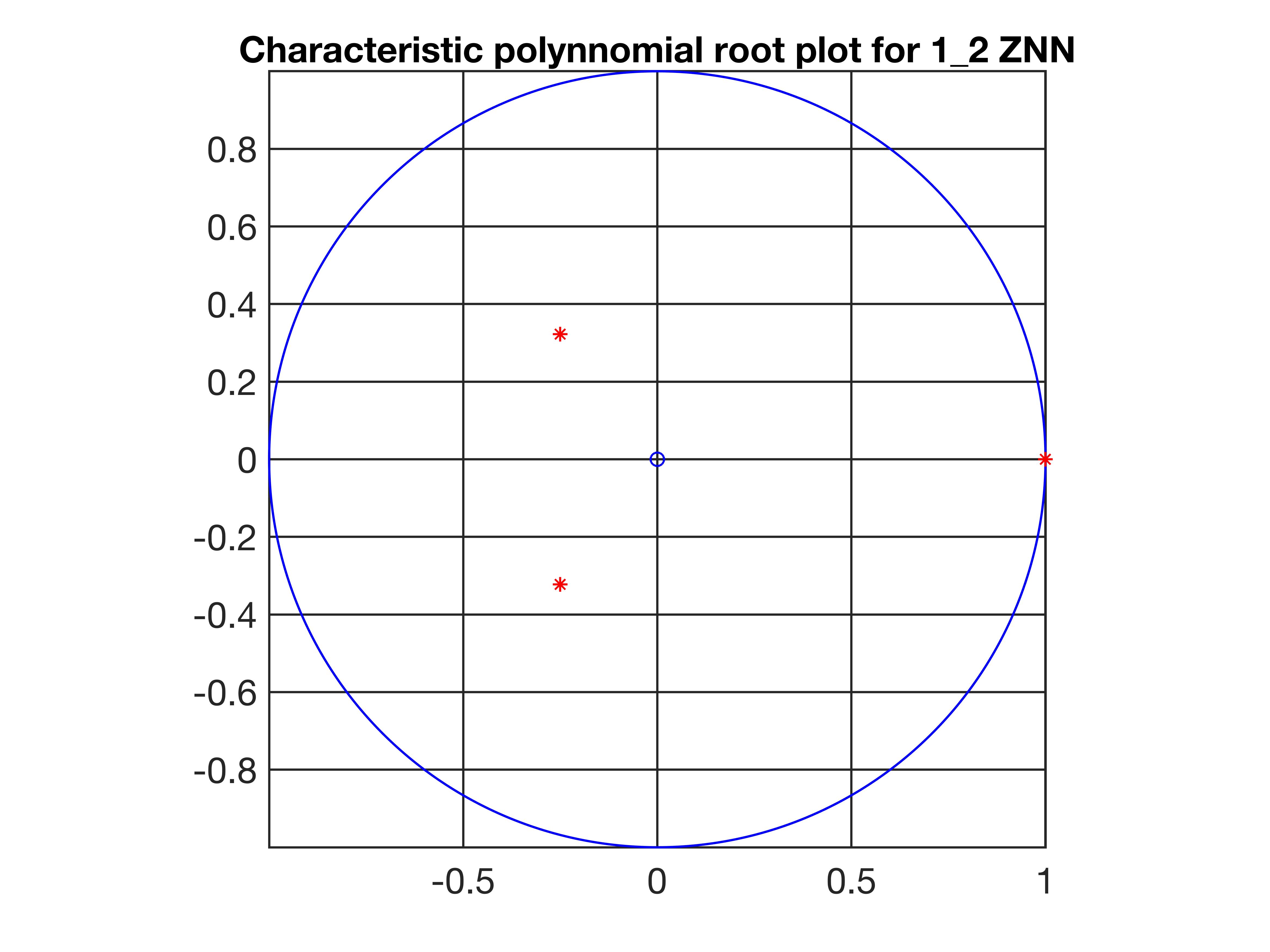} \includegraphics[width=79mm]{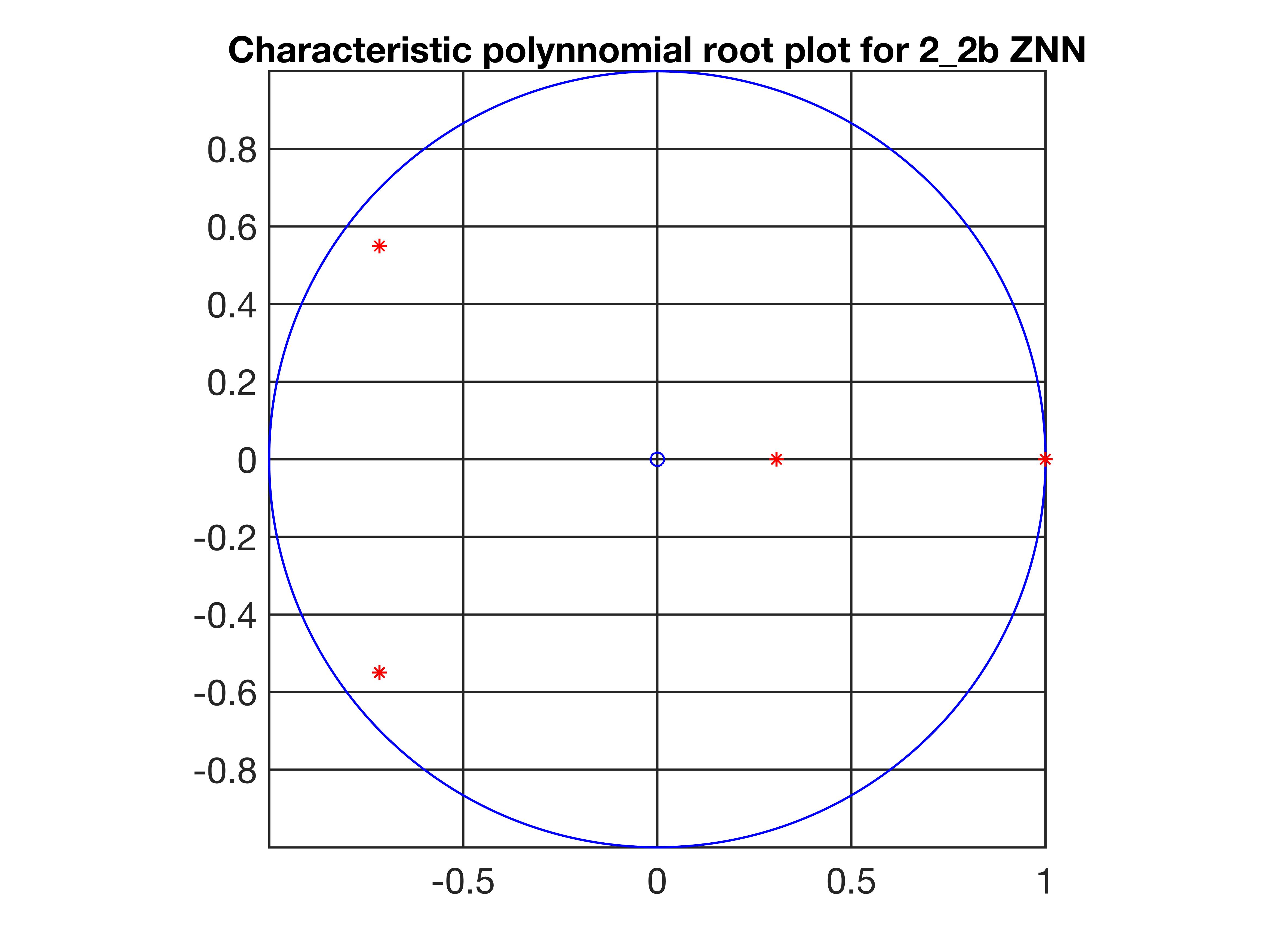}\\[-3mm]
\includegraphics[width=79mm]{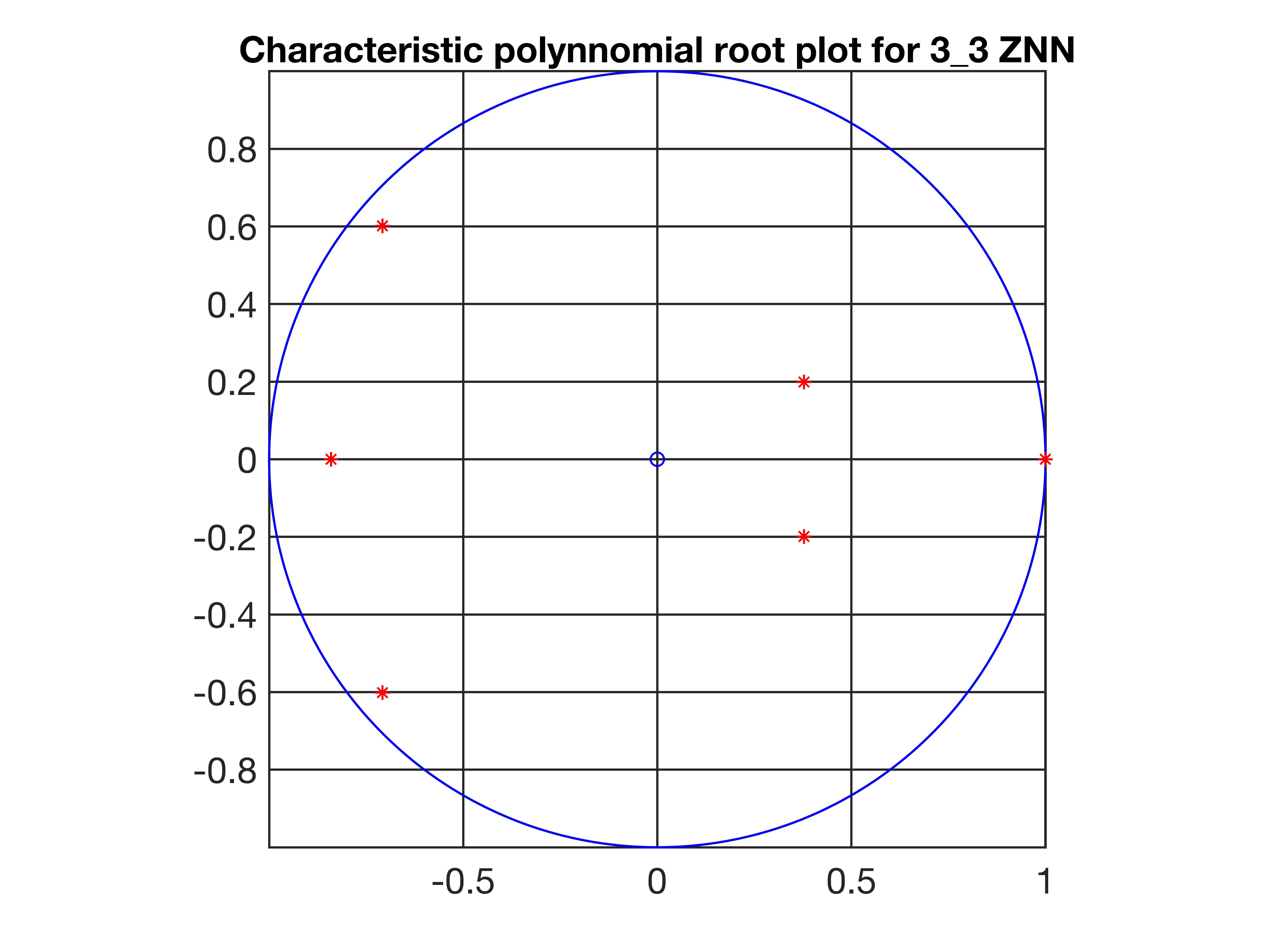} \includegraphics[width=79mm]{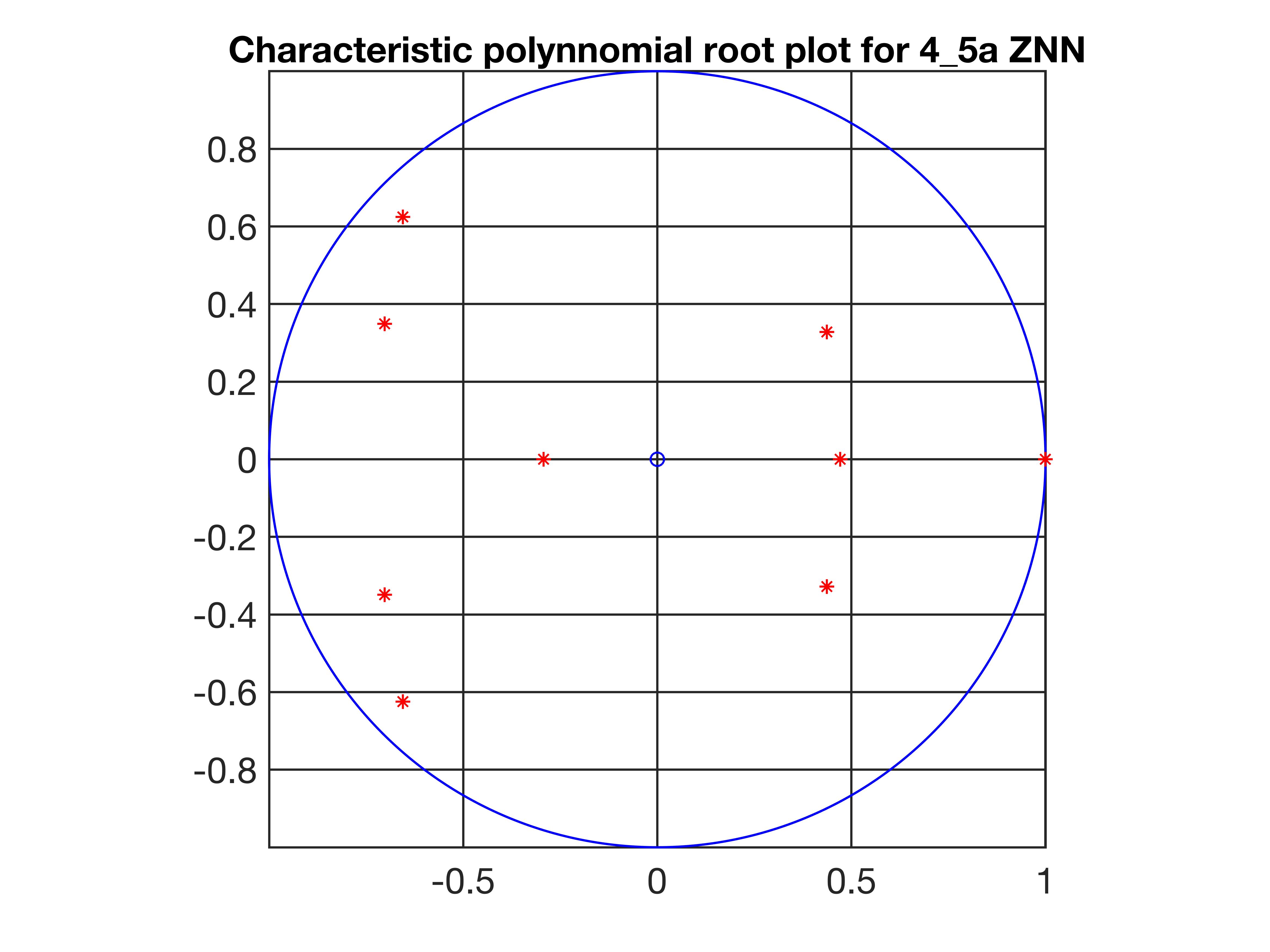}\\[-3mm]
\includegraphics[width=79mm]{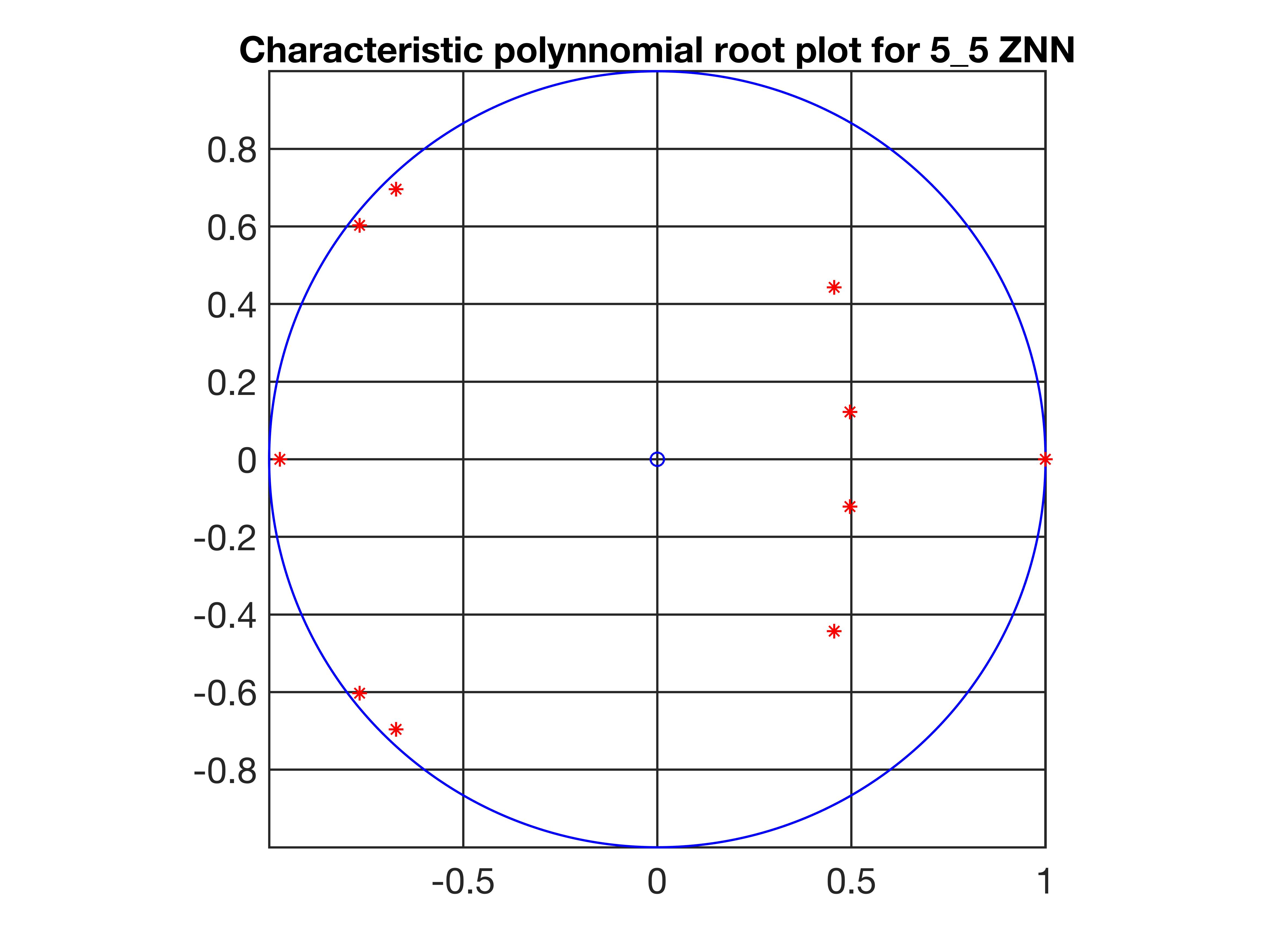} \includegraphics[width=79mm]{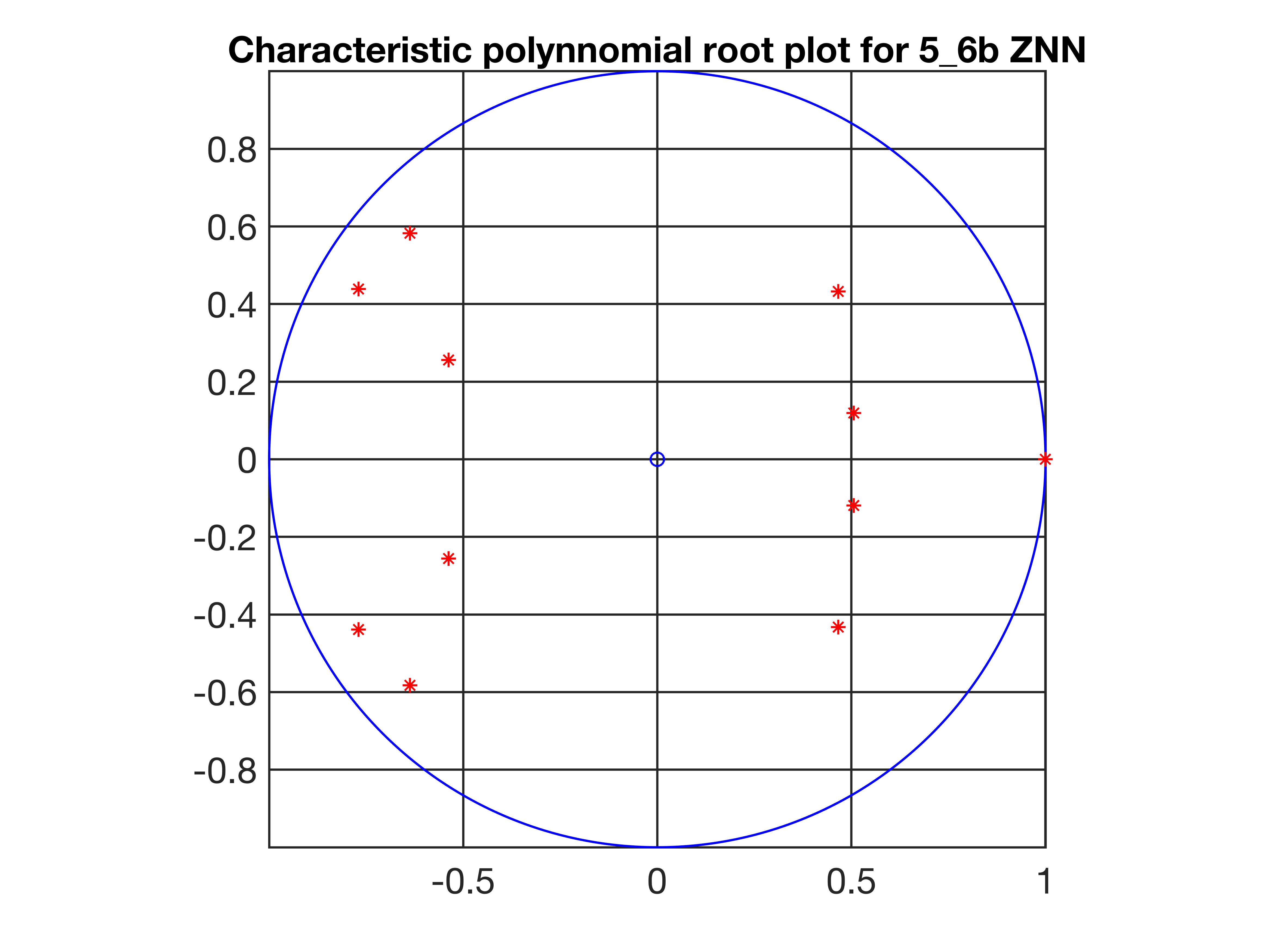}\\[-4mm]
Figure 2
\end{center}  
\normalsize

\noindent
The above root plots for the characteristic polynomials of \verb^5_5^ and  \verb^5_6b^ with identical truncation error order 7 show a slightly different geometry: five  of the roots  of \verb^5_5^ in the open unit disk lie much closer to the unit circle than those four inside roots closest to the unit circle for \verb^5_6b^. Moreover, the roots of \verb^5_6b^ are more widely and more evenly spread  over the unit disk than those of \verb^5_5^. In our computations, see Table 3, the results with ZNN and \verb^5_6b^ are much improved over those with ZNN and \verb^5_5^ for $\tau = 0.001$ and $\tau =0.0001$ with quite different values for $h$ for the identical random entry matrix $A_{243,243}$. Is there a correlation between a polynomial's root geometry, its accuracy, and the optimal value of $h$ in the associated finite difference ZNN scheme? If so, which geometries and which values of $h$ give better results; and why? These are open questions.\\[-3mm]

Our data was computed using MATLAB\_R2018a on a MacBook Pro (early 2015 version) under High Sierra OS with a 4 core Intel i5 processor and 16 GB of Ram.
The MATLAB codes {\tt FOVZNN4\_5ashorteigsym.m} for method \verb^4_5a^ and {\tt FOVZNN5\_6bshorteigsym.m} for method \verb^5_6b^ can be found and downloaded from \cite{FUFOVZNN45}.\\[2mm]
{\bf Remark :} \emph{( On repeated eigenvalues inside symmetric matrix flows $A(t)$)}
\\[0.8mm]
\hspace*{6mm} \begin{minipage}{153mm}{For both, decomposition methods and path following ODE initial value algorithms that find symmetric matrix flow eigendata over time, it is crucial that the matrix flow $A(t)$ never encounters repeated eigenvalues in its course, see Dieci and Eirola \cite{DE99}  and  Loisel and Maxwell \cite{LM}, respectively. 
 In \cite{LM} elaborate work-arounds for symmetric matrix eigenvalue repetitions are developed that can mitigate the built-in method limitations.\\
Our discretized ZNN eigenvalue  method encounters no problems when faced with matrix flows that incur repeated eigenvalues. If we had computed all eigenvalues of each $A(t_k)$ instead of just the two  extreme ones and then at each step had  selected the two extreme eigenvalues and used their associated eigenvectors to find the corresponding FoV points,  this extension of our algorithm and code would deal securely with such derogatory symmetric matrix flows.\hfill
That task sounds like a good project for a strong master's candidate.}
\end{minipage}\\[2mm]
{\bf Conclusions :}\\[0.8mm]
\hspace*{6mm} \begin{minipage}{153mm}{The ZNN method is completely new and highly accurate for time-varying matrix problems. It differs from all other time-varying matrix methods in its use of an error equation inspired differential equation that assures ZNN of automatic exponentially fast convergence. This new ODE is then solved by forward looking finite difference equations  that have never before been used and were actually unknown before their appearance with ZNN. Being completely new and different from decomposition methods such as Dieci et al \cite {DE99} or path following ODE based methods such as \cite{CU92, LM, WCW16, YFT04},  ZNN bypasses some of the erstwhile algorithm caused restrictions and limitations. It often works well in wider, less resticted  applications than was possible before. \\[1mm]
ZNN is a newly emerging  worthy subject of numerical analysis, ready made  for thorough studies that may answer its many  new open questions. Our community is thus called to set ZNN on firm computational and numerical ground.}
\end{minipage}
\newpage

\noindent
{\bf  Afterthoughts :} \\[1mm]
\hspace*{6mm}\begin{minipage}[t]{153mm}{When I was younger, much younger, I sometimes envisioned of someday replacing how we computed matrix eigenvalues and eigenvectors then, namely via matrix factorizations such as QR or via vector and subspace iteration schemes, with solving large numbers of linear equations instead. That was at times when I was trying to gain deeper understandings; guessing at what might lie ahead for us in  'scientific computing' terms as we would call it now and well before that term was invented. I mused about it, but found no answers then.\\[1mm]
This paper's main algorithm  spends most of its time (up to nearly half of its total CPU time)  on solving thousands and thousands of linear equations in its most expensive MATLAB code line \verb^X = linsolve(P,q)^. On the other hand,  its finite difference scheme formulated  eigendata finding part  runs almost effortlessly and in almost no time with   linearly combining vectors cleverly at $O(n)$ cost. As applied to the FoV boundary problem here, Zhang Neural Networks  can find time-varying and parameter-varying eigendata in record time and record accuracy today by using cheap and  clever linear vector combinations while relying on still costly  linear equation solves instead of Francis multishift QR.  \ Amazing.\\[1mm]
Does this answer my day-dreams of old, I wonder. \ What else is in store for us in future matrix numerics?\\[1mm]
And : \emph{Google} no longer uses the SVD and Francis QR,  as Gene Golub had  once  suggested in Google's `garage days', nor does \emph{Google}  use its 'page rank algorithm' as search engine any longer. Instead it uses  a different  Neural Network based search engine, called TensorFlow, since around 2013, see e. g. \cite{S2018}.}
\end{minipage}

\vspace*{10mm}

\noindent
\centerline{{[} .. /latex/FOVviaZNN.tex] \quad \today }

\vspace*{9mm}
\begin{minipage}{6cm}
{compareLM\_ZNN.jpg\\
roots12.jpg\\
roots22b.jpg\\
roots33.jpg\\
roots45a.jpg\\
roots55.jpg\\
roots56b.jpg\\[2mm]
7 image files}
\end{minipage}

\end{document}